\documentclass[12pt]{amsart}
\usepackage{amssymb,latexsym}
\input epsf
\usepackage{graphicx}
\newtheorem{theorem}{Theorem}[section]
\newtheorem{proposition}[theorem]{Proposition}
\newtheorem{corollary}[theorem]{Corollary}
\newtheorem{definition}[theorem]{Definition}

\newtheorem{lemma}[theorem]{Lemma}
\newtheorem{example}[theorem]{Example}

\newcommand{\reals}{\mathbb{R}}
\newcommand{\inv}{\mathrm{Inv}}
\newcommand{\invr}{\mathrm{Inv}_{\rho}}

\begin{document}
\title[Climbing elements]
{Climbing elements in finite Coxeter groups}

\author[Brady]{Thomas~Brady}
\address{School of Mathematical Sciences\\
Dublin City University\\
Glasnevin, Dublin 9\\
Ireland}
\email{tom.brady@dcu.ie}

\author[Kenny]{Aisling~Kenny}
\address{School of Mathematical Sciences\\
Dublin City University\\
Glasnevin, Dublin 9\\
Ireland}
\email{aisling.kenny9@mail.dcu.ie}

\author[Watt]{Colum~Watt}
\address{School of Mathematical Sciences\\
Dublin Institute of Technology\\
Dublin 8\\
Ireland}
\email{colum.watt@dit.ie}

\date{}

\begin{abstract}
We define the notion of a climbing element in a finite real
reflection group relative to a total order on the reflection set and we characterise
these elements in the case where the total order 
arises from a bipartite Coxeter element.
\end{abstract}

\maketitle
\section{Introduction}
Suppose $(W, S)$ is a finite Coxeter system.   Each reduced expression for an element
$w$ of $W$ determines a total order on the inversion set of $w$.  
The inversion set of the  longest element $w_0$ of $W$ is equal to 
the set, $T$,  of all the reflections 
and a particular reduced expression for $w_0$ gives a total order, $\le_T$, on  $T$.
For some elements $w$ of $W$, the restriction of $\le_T$ to the inversion set of $w$
coincides with the order determined by one of its reduced expressions.
We will call such an element $w$ a \emph{climbing} element of $W$.
Geometrically, this means that there is a gallery from the fundamental domain $C$ to $w(C)$
which crosses hyperplanes in increasing order.
\vskip .2cm
In this paper, we characterise the climbing elements in the case where the reduced
expression for $w_0$ is obtained by iterating a bipartite factorisation of a Coxeter element.
This characterisation is obtained using the construction from \cite{BW-P2A} of a copy
of the type-$W$ generalised associahedron, realised as a coarsening of the fan
determined by the $W$ reflection hyperplanes.   This coarsening determines an
equivalence relation on $W$ whose equivalence classes we prove directly to be
intervals in the left  weak order. The least elements of these intervals
are precisely the climbing elements. 
It follows that the number of climbing elements
is equal to the $W$-Catalan number.   
The maximal elements in these
intervals are translates of the falling elements of $W$,
a notion that is analagous to that of climbing elements
but which is defined using
the reverse of the order $\le_T$.
\vskip .2cm
For any minimal factorisation of a Coxeter element,  the interval property of the equivalence classes can be deduced from \cite{RedCamb} and \cite{RS1}, where the minimal elements are the corresponding Coxeter-sortable elements.  Thus we provide a different characterisation of Coxeter-sortable elements in the case of a bipartite factorisation of the Coxeter element. 
\vskip .2cm
The paper is organised as follows.  In \S $2$ we collect some facts about inversion sets, extend a theorem of Papi and
recall some results from \cite{BW} and \cite{BW-P2A} about orderings of roots and the geometry of
the generalised associahedron.
In \S $3$ we define climbing elements and we show that each facet of  the generalised associahedron
determines such an element.  We characterise climbing elements  in \S4 while in \S5 we introduce and characterise
falling elements.
\section{Preliminaries}
\subsection{Inversion sets}
\label{Inv}
 For background on reflection groups, root systems
and inversion sets we refer to \cite{BB} and \cite{Bou}.
Throughout this paper, $(W,S)$ is a Coxeter system  with $W$ finite, acting effectively on $\reals^n$
and with standard generating set $S = \{s_1, \dots s_n\}$. Denote by $T$ the reflection set
of $W$, that is, the set of congugates of elements of $S$.  Let $C$ be the fundamental
chamber with inward unit normals given by the simple roots 
$\{\alpha_1, \dots , \alpha_n\}$, 
where $s_i$ is the reflection in the hyperplane normal to $\alpha_i$.  
Let  $\{\beta_1, \dots , \beta_n\}$ be the dual basis so that 
$\alpha_i\cdot \beta_j = \delta_{ij}$.  
\vskip .2cm
For each $w \in W$ we define $\invr(w)$ to be the set of positive roots $\lambda$ such that
$w^{-1}(\lambda)$ is a negative root.  Thus
$\invr(w)$ is the set of positive roots whose orthogonal 
hyperplanes separate the fundamental chamber $C$
from its image $w(C)$.  The corresponding set of reflections is denoted by $\inv(w)$,  that is,
$\inv(w) = \{R(\lambda) \mid \lambda \in \invr(w)\}$, 
where $R(\lambda)$ is the reflection in the hyperplane orthogonal to $\lambda$.   We refer to $\inv(w)$ as the \emph{inversion set} of $w$.
If
$w = s_{i_1}s_{i_2} \dots s_{i_k}$  is a  reduced word, then $\inv(w)  = \{t_1, \dots , t_k \}$ where
\begin{equation}
\label{D-t_i}
t_1 = s_{i_1},\ \ t_2 = s_{i_1}s_{i_2}s_{i_1},  \ \ t_3 = s_{i_1}s_{i_2}s_{i_3}s_{i_2}s_{i_1},\ \ \dots \ \ 
\end{equation}
as in section 1.3 of \cite{BB}.   This defines a linear order on $\inv(w)$ and the corresponding linear order on $\invr(w)$ is given by
\[\alpha_{i_1},\  s_{i_1}(\alpha_{i_2}),\  s_{i_1}s_{i_2}(\alpha_{i_3}) , \ \dots \ , \ s_{i_1}\dots s_{i_{k-1}}(\alpha_{i_k}).\]
In \cite{P}, Papi characterises ordered inversion sets among ordered subsets of
$T$.  His proof is given for crystallographic groups although he notes that it can be
generalised to apply to all Coxeter groups. The following theorem modifies Papi's
characterisation and is valid in the general finite case. For
completeness the proof is included  in an appendix.
\begin{theorem}
An ordered subset $\Sigma$ of positive roots for $W$ is derived from a reduced expression for some element of
$W$ if and only if $\Sigma$ satisfies both of the following conditions on triples $\{\sigma, \tau, \rho\}$ of positive roots 
satisfying $\rho = a \sigma+ b \tau$ for some numbers $a > 0$ and $b > 0$.\\
(i)  Whenever $\sigma$ and  $\tau$ are elements of $\Sigma$ with $\sigma < \tau$ then $\rho \in \Sigma$ and $\sigma < \rho < \tau$.
\\
(ii)  Whenever $\rho$ is an element of $\Sigma$, then either {\it (a)} $\sigma \in \Sigma$ and $\sigma < \rho$
or {\it (b)} $\tau \in \Sigma$ and $\tau < \rho$.
\label{Papi}
\end{theorem}
\subsection{Geometry of the generalised associahedron}
\label{assocgeom}
Let $\ \Pi^+$ be the set of all positive roots.   
We recall from \cite{BW} the special features
of the linear order on $\ \Pi^+$ determined by iterating a bipartite Coxeter element.
First assume that the elements of the simple system are
ordered so that $\{\alpha_1, \dots , \alpha_s\}$
and $ \{\alpha_{s+1}, \dots , \alpha_n\}$ are orthonormal
sets. Let $c = R(\alpha_1) R(\alpha_{2}) \dots R(\alpha_n)$ be the corresponding
bipartite Coxeter element and let $h$ denote the order of $c$.    If $w_0$ is the longest element
of $W$ then it follows from the proof of Corollary~4.5 of \cite{S} that $w_0$ has the  reduced expression
\[w_0 = \left\{\begin{array}{cr}
c^{h/2} & \mbox{if $h$ is even\ }\\
c^{(h-1)/2}R(\alpha_1)  \dots R(\alpha_s) & \mbox{if $h$ is odd.}
\end{array}\right.
\]
It follows that the ordered set $\inv_\rho(w_0)$ is equal to  $\{\rho_1, \rho_2, \dots , \rho_{nh/2}\}$ where
\[ \rho_i = R(\alpha_1)R(\alpha_2)\dots R(\alpha_{i-1})\alpha_i,\]
and the $\alpha$'s are indexed cyclically modulo $n$.   In fact, $\inv_\rho(w_0) = \Pi^+$ and we denote
this order on $\Pi^+$ by $\le_\rho$ and by $\le_T$ the corresponding order on the reflection set $T$.
\vskip .2cm
Furthermore,  in  \cite{BW},  we define the vectors
\[ \mu_i = R(\alpha_1)R(\alpha_2)\dots R(\alpha_{i-1})\beta_i
, \ \  i = 1, 2, \dots , nh\]
where the $\alpha$'s and the $\beta$'s are indexed cyclically modulo $n$.   It is immediate from the
definitions of $\rho_i$ and $\mu_j$ that
$\rho_{i+n} = c(\rho_i)$, $\mu_{j+n} = c(\mu_j)$ and $\rho_i \cdot \mu_i = 1$.  We recall that
$\mu_i = \mu(\rho_i)$ where $\mu$ is the linear map defined by $\mu = 2(I-c)^{-1}$.    In particular,
$\rho_i = (1/2)(I-c)\mu_i$.
Furthermore, we have
\vskip .2cm
\begin{proposition} \rm{(Proposition~4.6 of \cite{BW})}\\
(a)  $\mu_i\cdot \rho_{j} = -\mu_{j+n} \cdot \rho_{i}$ for all $i$ and $j$.\\
(b) $\mu_i\cdot \rho_j \ge 0$, for $1 \le i \le j \le nh/2$.\\
(c) $\mu_{i+t}\cdot \rho_i = 0$, for $1 \le t \le n-1$ and
for all $i$.\\
(d) $\mu_{j} \cdot \rho_i \leq 0$ for $1 \le i < j \le nh/2$.
\label{table}
\end{proposition}
We recall from \cite{BW-P2A} that a copy, denoted $\mu AX(c)$, of the type-$W$ associahedron is identified with
a particular coarsening of the Coxeter fan, that is, of the fan defined by the $W$ reflection hyperplanes.
This coarsening has rays in the directions $\mu_1, \dots , \mu_{nh/2+n}$
and a facet with vertex set $\{\mu(\tau_1), \dots , \mu(\tau_n)\}$
whenever both
\[\rho_1 \le \tau_1 < \tau_2 < \dots < \tau_n \le \rho_{nh/2+n} \ \ \ \mbox{and}\ \ \ c = R(\tau_n)\dots R(\tau_1).\]
\vskip .2cm
We define an equivalence
relation on $W$ by $w \sim w'$ if and only if $w(C)$ and $w'(C)$ are
contained in the same facet of $\mu AX(c)$.
\vskip .2cm
Finally, we will use the filtration of $\mu AX(c)$ inherited from the filtration of $X(c)$ used in \cite{BW}.  For
each root  $\rho$ we define the subsets $\rho^+$, $\rho^-$ and $\rho^{\perp}$ by
\begin{eqnarray*}
\rho^+ &=& \{ x \in \reals \mid  x \cdot \rho \ge 0 \}\\
\rho^- &=& \{ x \in \reals \mid  x \cdot \rho \le 0 \}\\
\rho^{\perp} &=& \{ x \in \reals \mid  x \cdot \rho =  0 \}.
\end{eqnarray*}
For $n \le i \le nh/2+ n$, we define $V_i = \{\mu_1, \dots , \mu_i\}$, $\mu X_i$ to be  the subcomplex of $\mu AX(c)$
consisting of those simplices with vertices in $V_i$ and
\[
\mu Z_i = \rho_{i-n+1}^+ \cap \rho_{i-n+2}^+ \cap \dots \cap \rho_{nh/2}^+.
\]
It follows that the closure of $\mu Z_i \setminus \mu Z_{i-1}$  is equal to
\[\rho_{i-n}^{-}\cap \rho_{i-n+1}^{+}\cap \dots \cap \rho_{nh/2}^{+}.\]
We note that $\mu Z_n$ and $\mu Z_{nh/2+n}$ coincide with the fundamental chamber  $C$  and with $\reals^n$ respectively.  
We also note that Proposition~7.6 of \cite{BW} (in the case $\alpha = c$) can be extended to show that 
$\mu Z_i $ coincides with both the positive cone on $\mu X_i$ and the positive span of $V_i $.
\section{Climbing elements}
\label{secclimb}
In this section we define climbing elements and show that 
each subset of the vertex set of $\mu AX(c)$
determines a climbing element.  In the case of the vertex set of a facet we will show that this climbing element is 
the minimum in the corresponding equivalence class of $(W, \sim)$.
\begin{definition}
An element $w$ of $W$ is climbing if the order on $\mbox{Inv}(w)$ given by $\le_T$
coincides with the order determined by one of the reduced expressions for $w$.
\label{climb}
\end{definition}

\begin{definition}
For each subset $A$ of $V_{nh/2+n}$ we define the set $N(A)$ of positive roots by
\[ \mbox{N}(A) = \{\rho_i \mid 1 \le i \le nh/2 \mbox{ and }
 \rho_i \cdot \mu \le 0 \mbox{ for all } \mu \in A\}.\]
\label{partition}
\end{definition}
Thus a positive root $\rho$ belongs to $\mbox{N}(A)$  if and only if 
$A \subseteq \rho^-$.
\begin{example}  If $A = \{\mu(\rho_i) \}$ then Proposition~\ref{table} implies that
$$\mbox{N}(A)=\{\rho_j : j < i \mbox{ or } \rho_j \cdot \mu(\rho_i) = 0 \}.$$ 
For a larger set $B$, $\mbox{N}(B)$ is the intersection of sets of this form.
\label{membership}
\end{example}
\begin{proposition}  For each subset $A$ of 
$V_{nh/2+n}$ there exists a element $w \in W$  such that the ordered set $(\mbox{N}(A), \le_\rho)$ coincides with
the ordered set $\invr(w)$ for some reduced expression of $w$.    In particular, $w$ is climbing.
\label{genclimb}
\end{proposition}
\textbf{Proof:} We show that $N(A)$ satisfies the criteria (i) and (ii) of Theorem~\ref{Papi}. First
suppose $\rho_i, \rho_j \in \mbox{N}(A)$ with $i < j$ and that $a, b > 0$ are such that
$\rho_k = a\rho_i + b\rho_j$ is a positive root.
For each $\mu \in A$ we have
\[\rho_k \cdot \mu = (a\rho_i + b\rho_j)\cdot \mu =
a(\rho_i\cdot \mu) + b(\rho_j\cdot \mu) \le 0\] since
$\rho_i, \rho_j \in N(A)$.   Thus, $\rho_k \in \mbox{N}(A)$.   As the order $\le_\rho$ on $\Pi^+$ is derived from
a particular reduced expression for the longest element $w_0$, the `only if' part of Theorem~\ref{Papi} yields
$\rho_i \le_\rho \rho_k \le_\rho \rho_j$ and criterion (i) follows.

Next, suppose that $\rho_i$ and $\rho_j$
are positive roots with  $i < j$ and that $a, b > 0$ are such that
$\rho_k = a\rho_i + b\rho_j \in \mbox{N}(A)$.
As in the previous paragraph, Theorem~\ref{Papi}  yields
 $\rho_i \le_\rho \rho_k \le_\rho \rho_j$.   It remains to
show that $\rho_i \in \mbox{N}(A)$.

If $\rho_i \not\in \mbox{N}(A)$ then $\rho_i \cdot \mu > 0$
for some $\mu \in A$.  By definition of $\mu AX(c)$, $\mu = \mu(\rho_q)$ 
for some root $\rho_q$ with $1 \le q \le nh/2+n$.
In fact, $1 \le q \le nh/2$ since $\{\mu_{nh/2+1}, \dots , \mu_{nh/2+n}\}$ 
are the rays of the cone $w_0(C)$, the opposite
chamber to $C$.  Now part (d) of Proposition~\ref{table} gives $q \le i$.   
Therefore $q < j$ and, hence, part (b) of Proposition~\ref{table}
implies that $\rho_j \cdot \mu  = \rho_j \cdot \mu(\rho_q) \ge 0$.   Thus
\[\rho_k \cdot \mu = (a \rho_i + b \rho_j)\cdot \mu  = 
a (\rho_i \cdot \mu) + b (\rho_j\cdot \mu) \ge a (\rho_i \cdot \mu) > 0,\]
contradicting the assumption that $\rho_k \in \mbox{N}(A)$.
\vskip .2cm
If $F$ is a facet of $\mu AX(c)$, 
we denote its set of vertices by $V_F$. That is  $V_F= F \cap V_{nh/2+n}$. 
Such vertex sets will be particularly
important in the sequel. 
 \begin{proposition}
If  $F$ is a facet of $\mu AX(c)$ and $x_F \in W$ is the climbing element with
$\invr (x_F) = N(V_F)$, then $x_F(C) \subseteq F$.
\label{containment}
\end{proposition}
\textbf{Proof:}    Since $\mu AX(c)$ is a coarsening of the Coxeter fan, 
the facet~$F$ can be characterised as an intersection
of halfspaces determined by the roots $\rho_i$.  We show that $x_F(C)$ 
is contained in the same intersection.   If  $\rho_i$ is a positive root
with $F$ contained in $\rho_i^-$, 
then $\rho_i \in N(V_F)$.  Since $N(V_F) = \invr(x_F)$,  it follows that 
$x_F(C)$ must also be contained in $\rho_i^-$.  
On the other hand, if $\rho_j$ is a positive root  with $F$  contained in
$\rho_j^{+}$ then $\rho_j \not\in N(V_F)$ since $F$ has nonempty interior and 
hence cannot be contained in $\rho_j^\perp$.  Thus $x_F(C)$ must also be contained in 
$\rho_j^{+}$.

\begin{corollary}
Each equivalence class of $(W,\sim)$ contains a minimum in the
left weak order on $W$.
\label{interval}
\end{corollary}
\textbf{Proof:}   Let $F$ be a facet of $\mu AX(c)$ with vertex set $V_F$
and let $x_F$ be the element of $W$ whose inversion set is $N(V_F)$ (Proposition~\ref{genclimb}).    
By Proposition~\ref{containment},  $x_F(C)$ is contained in $F$.
If $w \sim x_F$ then $w(C) \subseteq F$, by definition, and it
follows that 
$w(C) \subset \rho_i^-$ for each $\rho_i \in N(V_F)$. Thus  
$N(V_F) \subseteq \invr(w)$ and
Proposition~3.1.3 of \cite{BB} now implies that $x_F$ precedes $w$ 
in the left weak order on $W$.
\vskip .2cm

\section{Characterising climbing elements}
\label{char}
The proof of Corollary~\ref{interval} shows that
the number of  facets of $\mu AX(c)$ does not exceed the number of climbing elements.  
In fact the  theorem below implies that these two numbers are equal.

\begin{lemma}
\label{L41}
If $\mu(\rho_i)$ is the last vertex of a facet $F$  of $\mu AX(c)$ and
if $w$ is a climbing element  for which $w(C) \subset F$,
then $R( \rho_{i-n})w$ is also
a climbing element.
\end{lemma}
\textbf{Proof:} 
Assume that 
$\mu(\rho_{i_1})$, $\mu(\rho_{i_2})$, \dots, $\mu(\rho_{i_{n-1}})$, $\mu(\rho_{i})$
are the vertices of $F$ where
$1 \le i_1 <  \dots < i_{n-1}<i \le nh/2+n $ and
$c = R(\rho_{i})R(\rho_{i_{n-1}})\dots R(\rho_{i_1})$.
Since  $1 \le i-n \le nh/2$ and
\[ c = R(\rho_{i})R(\rho_{i_{n-1}})\dots R(\rho_{i_1}) = 
R(\rho_{i_{n-1}})\dots R(\rho_{i_1})R( \rho_{i-n}), \]
Lemma~2.2 of \cite{ABMW} implies that  $ \rho_{i-n} \cdot \mu(\rho_{i_k}) = 0$  
for $k = 1, 2, \dots , n-1$. 
Thus  the face  of $F$ opposite to the vertex $\mu(\rho_{i})$
is contained in the hyperplane $ \rho_{i-n} ^{\perp}$. 
It follows that 
\[ F \subseteq \overline{\mu Z_{i}\setminus \mu Z_{i-1}}=
\rho_{i-n}^{-}\cap \rho_{i-n+1}^{+}\cap \dots \cap \rho_{nh/2}^{+} \]
and, hence, 
the last wall crossed by any increasing gallery for $w$
is $ \rho_{i-n} ^{\perp}$.  If we delete the last chamber from such an increasing
gallery, we obtain an increasing gallery for $R( \rho_{i-n})w$.
Therefore $R( \rho_{i-n})w$ is a climbing element, as required.

\begin{theorem}
Each equivalence class of $(W,\sim)$ contains exactly one 
climbing element.  In particular the number of climbing elements 
is equal to the $W$-Catalan number.
\label{uniqueness}
\end{theorem}
\textbf{Proof:}    Fix an associahedron facet $F$ whose vertices
are $\mu(\rho_{i_1})$, $\mu(\rho_{i_2})$, \dots, $\mu(\rho_{i_{n-1}})$, $\mu(\rho_{i})$
where
$1 \le i_1 <  \dots < i_{n-1}<i \le nh/2+n $ and
$c = R(\rho_{i})R(\rho_{i_{n-1}})\dots R(\rho_{i_1})$.
We need to show that there is only one 
climbing element $w \in W$ for which $w(C) \subset F$.
Our proof is by 
induction on $i$.
\vskip .2cm
First note that $i \ge n$ and 
if $i = n$ then  $F$ must coincide with the
fundamental domain~$C$. In this case the
identity element of $W$ is the only element for which
$w(C) \subset F$. 
\vskip .2cm
Assume now  that $i > n$ 
and that for each associahedron facet $F'\subseteq \mu Z_{i-1}$
there is a unique climbing element $w'$ for which $w'(C) \subseteq F'$. 
Let $G$ be the only other associahedron facet which
contains  the face $F \cap \rho_{i-n} ^{\perp}$.
Since 
$\mu(\rho_i) \cdot  \rho_{i-n} <  0$,  $G$ is contained in $ \mu Z_{i-1}$.
Then $[R( \rho_{i-n})w](C)$  also lies in
$G$ since $G$ shares the face $ \rho_{i-n} ^{\perp} \cap F$ with
$F$. As $R( \rho_{i-n})w$ is climbing (by Lemma~\ref{L41}),
the induction hypothesis implies that
$R( \rho_{i-n})w  = w'$, the unique climbing element for which $w'(C) \subseteq G$.
Hence $w = R( \rho_{i-n}) w'$ is uniqely determined.
\begin{corollary}   The set of climbing elements in $W$ coincides with the set of 
Coxeter-sortable elements of $W$.
\label{coxsort}
\end{corollary}
\textbf{Proof:}   By Theorem~1.1 of \cite{RedCamb} the Coxeter-sortable elements of $W$ are precisely the  minima of the equivalence classes 
of $(W, \sim)$.  By  Theorem~\ref{uniqueness} and the proof of Corollary~\ref{interval}, the climbing elements are also the  minima of these equivalence classes.   
\section{Falling elements}
In this section we show that each equivalence class of $(W,\sim)$ contains a maximum in the
left weak order on $W$.    Just as a climbing element is reached from the fundamental chamber $C$ via a gallery which crosses hyperplanes in 
increasing order, each of these maxima is reached from the opposite chamber $w_0(C)$ via a gallery which crosses hyperplanes in 
decreasing order.   In order to use the results of sections~\ref{secclimb} and \ref{char} our strategy is to 
rebuild  $\mu AX(c)$ with $w_0(C)$ taking the place of $C$ and $c^{-1}$ taking the place of $c$.  This will give an 
ordering on $T$ which is the reverse of the order $\le_T$ and we will refer to the corresponding notion of climbing element  
as a falling element.  The required maxima will then have the form $fw_0$ where $f$ is falling.  
\vskip .2cm
Since the inward pointing normals for $w_0(C)$ are just the negatives of the inward pointing normals for $C$, the new simple system will 
be $\{-\alpha_1, \dots , -\alpha_n\}$.   We will order this simple system by using the corresponding order on the dual basis.   Sometimes this 
order is different than  the order $-\alpha_n, \dots , -\alpha_1$ but we will see that it gives the reverse of the order $\le_T $ on $T$.
\begin{definition}
For $1\le j \le n$ we define $\beta'_i = \mu_{nh/2+n-i+1}$ and we define $\{\alpha'_1,  \dots , \alpha'_n\}$ to be the dual 
basis to $\{\beta'_1,  \dots , \beta'_n\}$.
\label{newdual}
\end{definition}
\begin{proposition}   $\{\beta'_1, \dots , \beta'_{n-s}\}$ is a permutation of $\{-\beta_{s+1}, \dots , -\beta_n\}$ and 
$\{\beta'_{n-s+1}, \dots , \beta'_{n}\}$ is a permutation of $\{-\beta_{1}, \dots , -\beta_s\}$.
\label{perm}
\end{proposition}
\textbf{Proof: }  This follows from Steinberg's  proof of Theorem~4.2 of \cite{S}, where the vectors he denotes by $\sigma$ and 
$\tau$ lie in the non-negative linear spans of our $\beta_1, \dots ,\beta_s$ and  $\beta_{s+1}, \dots ,\beta_n$, respectively.
\begin{corollary}   $\{\alpha'_1,  \dots , \alpha'_{n-s}\}$ is a permutation of $\{-\alpha_{s+1}, \dots , -\alpha_{n}\}$ while 
$\{\alpha'_{n-s+1},  \dots , \alpha'_{n}\}$ is a permutation of $\{-\alpha_{1}, \dots , -\alpha_{s}\}$.   In particular, 
$ c^{-1} = R(\alpha'_1)R(\alpha'_2)\dots R(\alpha'_n)$ is a bipartite factorisation.  
\label{cinverse}
\end{corollary}
\begin{definition}
With the convention that $\alpha'_{i+n} = \alpha'_i$ and $\beta'_{i+n} = \beta'_i$, we define $\mu'_i = R(\alpha'_1)R(\alpha'_2)\dots R(\alpha'_{i-1})\beta'_i$ and  $\rho'_i = R(\alpha'_1)R(\alpha'_2)\dots R(\alpha'_{i-1})\alpha'_i$.

\label{primes}
\end{definition}
Note that $\mu'_j = \beta'_j$ for $1 \le j \le n$ and $\mu'_{i+ n} = c^{-1}\mu'_i$.
\begin{proposition}  The vectors $\mu'_j$ and $\rho'_i$ are related to  $\mu_j$ and $\rho_i$ by
$\mu'_j = \mu_{nh/2+n-j+1}$ and $\rho'_i =  -\rho_{nh/2-i+1}$, 
for $1 \le j \le nh/2+n$ and $1 \le i \le nh/2$ respectively.
\label{correspond}
\end{proposition}
\textbf{Proof: }   For the first identity, write $j = mn+k$ with $0\le k < n$.  Then
\begin{eqnarray*}
\mu'_j &=& R(\alpha'_1)\dots R(\alpha'_{j-1})\beta'_j\\
&=& [c^{-1}]^mR(\alpha'_1)\dots R(\alpha'_{k-1})\beta'_k\\
&=& c^{-m}\beta'_k \ \ \mbox{since} \ \ \beta'_k \perp \alpha'_1, \dots, \alpha'_{k-1}\\  
&=& c^{-m}\mu_{nh/2+n-k+1}\\
&=& \mu_{nh/2+n-mn-k+1}\\
&=& \mu_{nh/2+n-j+1}.\\
\end{eqnarray*}
For the second identity, we use the first identity and the relationships 
\[c\mu_i = \mu_i -2\rho_i , \ \ \ c^{-1}\mu'_i = \mu'_i -2\rho'_i \]
to get 
\begin{eqnarray*}
2\rho'_i &=& (I-c^{-1})\mu'_i\\
&=& (I-c^{-1})\mu_{nh/2+n-i+1}\\
&=& (I-c^{-1})c\mu_{nh/2-i+1}\\
&=& -(I-c)\mu_{nh/2-i+1}\\
&=& -2\rho_{nh/2-i+1},
\end{eqnarray*}
for $1 \le i \le nh/2$.
\vskip .2cm
We now construct a copy of the type-$W$ associahedron using $c^{-1}$ instead of $c$ and $\{\mu'_1, \dots , \mu'_{nh/2+n}\}$ instead of 
$\{\mu_1, \dots , \mu_{nh/2+n}\}$.  We find that the geometric complex is exactly the same since the vertex sets coincide by Proposition
~\ref{correspond} and there is a facet on a set 
\[\{\mu'_{i_1}, \dots , \mu'_{i_n}\}\]
if and only if there is facet on the corresponding set  
\[\{\mu_{nh/2+n-i_1+1}, \dots , \mu_{nh/2+n-i_n+1}\}.\]
The reflection ordering $\rho'_1, \rho'_2, \rho'_3, \dots $ is the reverse of $\le_T$ and determines a different notion of climbing element which
we will now call  falling.
\vskip .2cm
\begin{definition}
An element $w$ of $W$ is falling if the order on $\mbox{Inv}(w)$ given by the reverse of the
total order $\le_T$ coincides with the order determined by one of the reduced expressions for $w$.
\label{falling}
\end{definition}
The results of sections~\ref{secclimb} and \ref{char} apply to give
\vskip .2cm
\begin{theorem}    Each equivalence class of $(W,\sim)$ determines a unique falling element $f$.
The element $fw_0$ is  the  maximal element in the corresponding equivalence class in the left weak order.
\end{theorem}
\begin{corollary}    Each equivalence class of $(W,\sim)$ is an interval in the
left weak order on $W$. 
\end{corollary}

\appendix
\section{}
Before proving Theorem~\ref{Papi}, we prove some elementary facts.
\begin{lemma}
\label{LMN3}
If the positive root $\rho$ is not  simple  then
we can write $\rho =   a \sigma +b\tau $
for some real numbers $a,b>0$ and some positive roots $\sigma, \tau$.
\end{lemma}
\vskip .2cm
\textbf{Proof: }  First, if  $\rho$ is any  positive root then
we can write
\[\rho = a_1\alpha_1 + \dots + a_n \alpha_n
\mbox{ \ \ \ with  \ \ } a_i \ge 0 \mbox{ \ for \ } 1\le i \le n \]
and it follows that
\[0 < \rho \cdot \rho = \rho \cdot \left(\sum_i a_i \alpha_i\right ) =
\sum_{i} a_i (\rho \cdot  \alpha_i).\]
yielding $\rho \cdot \alpha_i > 0$ for some simple root $\alpha_i$.
\vskip .2cm
Now suppose that $\rho$ is a non-simple, positive root and that
$\alpha_i$ is a simple root with $\rho \cdot \alpha_i >0$ as above.
Since  $\rho$ is not  a simple root, it follows that
$s_i (\rho) = \sigma$ is a positive root.
However, $\sigma = \rho- b \alpha_i$, and hence $\rho= \sigma+ b\alpha_i$
where $b=2(\rho \cdot \alpha_i) >0$, as required.
\vskip .2cm
Define the vector $v_0$ by $v_0 = \beta_1+  \dots + \beta_n$ and note that
$v_0$ lies in the
interior of the fundamental chamber $C$ since $v_0 \cdot \alpha_i = 1$, for each $i$.
Note also that for each $w \in W$, the set 
$\invr(w)$ is equal to  the set of positive roots $\lambda$ such that
$\lambda\cdot w(v_0) < 0$.
\vskip .2cm
\begin{lemma}
\label{LMneg}
If $w \in W$ and $w(\sigma) \in \invr(w)$ then $\sigma$ is a negative root.
\end{lemma}
\vskip .2cm
\textbf{Proof: }   Directly from the definition of $\invr(w)$ we have
\[\sigma \cdot v_0 = w(\sigma)\cdot w(v_0) < 0.\]
\vskip .2cm
\textbf{Proof of Theorem~\ref{Papi}: }  First, assume that the ordered set $\Sigma$  is derived from a reduced expression
$w = s_{i_1}s_{i_2} \dots s_{i_k}$ for some element $w \in W$.
Extend this to a reduced expression
\[s_{i_1}s_{i_2} \dots s_{i_{nh/2}}\]
for the longest element  of $W$, as in
Section~1.8 of \cite{Hu}.  For each $1 \le j \le nh/2$,   let
$w_j = s_{i_1}s_{i_2} \dots s_{i_j}$
be the $j$th prefix of this expression and note that $w = w_k$.

For condition (i),
assume that $\sigma < \tau$ are elements of $\Sigma$ and
that $\rho = a \sigma+ b \tau$ is a positive root for some $a,b>0$.
Then $R(\sigma) = t_{r_1}$, $R(\rho) = t_{r_2}$ and $R(\tau) = t_{r_3}$
for some $1\le r_1<r_3 \le k$ (by our assumption on $\Sigma$)
and some $1 \le r_2 \le nh/2$, and
where the  $t_j$ are given by equation~(\ref{D-t_i}).
We show that $r_1 < r_2 < r_3$ by eliminating the
other possibilities.  It then follows that $\rho \in \Sigma$.
First, if $r_2 < r_1 < r_3$ then
\[w_{r_2}(v_0)\cdot \sigma > 0 \ \ \mbox{and} \ \ w_{r_2}(v_0)\cdot \tau > 0 \ \
\mbox{while} \ \ w_{r_2}(v_0)\cdot \rho < 0.\]
This is impossible since $\rho$ is a
positive linear combination of $\sigma$ and $\tau$.   Similarly, if $r_1<r_3<r_2$ then
\[w_{r_3}(v_0)\cdot \sigma < 0 \ \ \mbox{and} \ \ w_{r_3}(v_0)\cdot \tau < 0 \ \
\mbox{while} \ \ w_{r_3}(v_0)\cdot \rho > 0\]
 which is also impossible since $\rho$ is a
positive linear combination of $\sigma$ and $\tau$.
\vskip .2cm
For condition (ii),
assume that $\sigma$ and $\tau$ are positive roots and that $a,b>0$ are such that
$\rho = a \sigma+ b \tau$ is an element of $\Sigma$. Thus $R(\rho) =t_r$,
for some $1 \le  r \le k$, and hence
 $w_r(v_0)\cdot \rho < 0$.  Since $\rho$ is a positive linear combination
 of $\sigma$ and $\tau$, at least one of $w_r(v_0)\cdot \sigma$ and
 $w_r(v_0)\cdot \tau$ must be strictly negative.
 Thus,  either $R(\sigma) \in \inv(w_r)$ and hence
 $\sigma \le \rho$ or $R(\tau)\in \inv(w_r)$ and hence
 $\tau \le \rho$.  As $a,b>0$, we can exclude the possibilities of $\sigma=\rho$
 or $\tau = \rho$.
\vskip .2cm
For the converse, assume that $\Sigma$ is a set of positive roots which
satisfies conditions (i) and (ii).  As in  \cite{P}, we proceed by induction on
the cardinality  of $\Sigma$.
To start the induction
we assume that $\Sigma=\{\rho\}$. It suffices to show
that $\rho$ is a simple root, for then
$(R(\rho)$ is the required group element.
If  $\rho$ is not  a simple root, then Lemma~\ref{LMN3}
implies that  $\rho=a\sigma+b\tau$ for some other positive roots $\sigma$ and $\tau$
and some $a, b>0$.
By condition (ii), either $\sigma$ or $\tau$ is also in  $\Sigma$,
contradicting the assumption that $\Sigma$ has cardinality one.
\vskip .2cm
For the inductive step, assume that $k>1$ and that
the result is true for sets of cardinality less than $k$.
 Assume that $\Sigma = \{\rho_1, \rho_2, \dots , \rho_k\}$
 satisfies conditions (i) and (ii). Then  the  ordered set
$\Sigma' = \{\rho_1, \rho_2, \dots , \rho_{k-1}\}$ 
also satisfies these two conditions and hence
there is a reduced expression $u = s_{i_1}s_{i_2} \dots s_{i_{k-1}}$,
of some element $u\in W$, such that
\[\rho_1 = \alpha_{i_1},  \ \   \rho_2 = s_{i_1}(\alpha_{i_2}), \ \ \dots , \ \ \rho_{k-1} = s_{i_1}s_{i_2} \dots s_{i_{k-2}}(\alpha_{i_{k-1}}).\]
If $u^{-1}(\rho_k)$ is a simple root, $\alpha_{i_k}$ say, then the positivity of
$ \alpha_{i_k}$ implies that
$l(s_{i_1}s_{i_2}\dots s_{i_k}) = l(u)+1$ and, hence,
$w=us_{i_k}=s_{i_1}s_{i_2}\dots s_{i_k}$ is the required minimal expression.
Thus it remains to show that
$u^{-1}(\rho_k)$ must be simple.
\vskip 0.2cm
Assume that $u^{-1}(\rho_k)$ is not simple.
As $\rho_k \not \in \invr(u)$, it follows that
$u^{-1}(\rho_k)$ is a positive root.
Then
$u^{-1}(\rho_k) = a\sigma +b\tau$
for some positive roots $\sigma$ and $\tau$ and some real numbers $a,b>0$, by Lemma~\ref{LMN3}.
Thus
\begin{equation}\label{E-PosLC}
\rho_k = au(\sigma) + bu(\tau) .
\end{equation}
In order to apply condition (ii) to this equation, we need to show that
neither  $u(\sigma) $ nor $u(\tau)$ can be a negative root.
For example, if $u(\sigma) $ is negative,  then
\[-u(\sigma)\cdot v_0 > 0 \ \ \ \mbox{and} \ \ \ -u(\sigma)\cdot u(v_0) = -\sigma\cdot v_0 < 0 \]
putting $-u(\sigma)$ in $\invr(u)$.  Thus $-u(\sigma) = \rho_i$ for some $i <k$.
Condition (i) applied to the  expression
$u(\tau) = (1/b)\rho_k+ (a/b)\rho_i$
now implies that
$u(\tau) \in \Sigma$ and $u(\tau) < \rho_k$,
so that, in fact,  $u(\tau) \in \Sigma' = \invr(u)$.   This gives a contradiction by Lemma~\ref{LMneg}.
\vskip .2cm
Thus both $u(\sigma) $ and $u(\tau)$ must be positive and by Lemma~\ref{LMneg} again neither belong to $\invr(u)$.
This gives a contradiction since condition (ii)
applied to equation~(\ref{E-PosLC}) implies that one of $u(\sigma) $ and $u(\tau)$  is in $\Sigma$ and precedes $\rho_k$,
putting one of $u(\sigma) $ and $u(\tau)$ in $\Sigma' = \invr(u)$.


\begin{thebibliography}{99}
%
\bibitem{ABMW}
C.A.~Athanasiadis, T.~Brady, J.~McCammond and C.~Watt,
\textit{$h$-vectors of generalized associahedra and noncrossing partitions},
Int. Math. Res. Not. 2006, Art. ID 69705, 28 pp
%
\bibitem{BB}
A.~Bj\"orner and F. Brenti,
\textit{Combinatorics of Coxeter groups},
Graduate Texts in Mathematics, vol. 231, Springer, New York, 2005.
%
\bibitem{Bou}
N.~Bourbaki,
\textit{Lie groups and Lie algebras. Chapters 7--9.} Translated from the 1975 and 1982 French originals by Andrew Pressley. Elements of Mathematics (Berlin). Springer-Verlag, Berlin, 2005.
%

\bibitem{BW}
T.~Brady and C.~Watt,
\textit{Noncrossing Partition Lattices in finite real reflection groups},
Trans. Amer. Math. Soc. 360 (2008), 1983-2005.
%
\bibitem{BW-P2A}
T.~Brady and C.~Watt,
\textit{From Permutahedron to Associahedron},
Proc. Edinburgh. Math. Soc. (to appear)
%
\bibitem{Hu}
J.E.~Humphreys,
Reflection groups and Coxeter groups,
Cambridge Studies in Advanced Mathematics {\bf~29},
Cambridge University Press, Cambridge, England, 1990.
%
\bibitem{P}
P.~Papi,
\textit{A characterization of a special ordering in a root system},
Proc. Amer. Math. Soc. 120 (1994), no. 3, 661--665.
%
\bibitem{RedCamb}
N.~Reading, 
\textit{Sortable elements and Cambrian lattices}
 Algebra Universalis 56 (2007), no. 3--4, 411--437. 

%
\bibitem{RS1}
N.~Reading and D.~Speyer,
\textit{Cambrian fans},
J. Eur. Math. Soc. (JEMS) 11 (2009), no. 2, 407--447.
%
\bibitem{S}
R.~Steinberg,
\textit{Finite reflection groups},
Trans. Amer. Math. Soc. 91, No. 3, (1959) 493-504.
\end{thebibliography}
\end{document}